\documentclass[reqno]{amsart}
\usepackage{amssymb,tikz}
\usepackage[T1]{fontenc}
\usepackage{enumitem,microtype,indentfirst,titlesec}
\setlist[enumerate]{format=\normalfont}
\usepackage{newtxtext,newtxmath}
\usepackage[justification=centering]{caption}
\usepackage[dvipsnames,svgnames]{xcolor}
\usepackage[pdflang=en-UK,colorlinks,urlcolor=Navy,linkcolor=Navy,citecolor=Navy]{hyperref}
\usepackage{upref}
\usepackage[backend=bibtex,bibstyle=alphabetic,citestyle=alphabetic]{biblatex}
\calclayout
\makeatletter
\def\section{\@startsection{section}{3}\z@{.8\linespacing\@plus.7\linespacing}{1\linespacing}{\normalfont\Large\bfseries \centering}}
\def\subsection{\@startsection{subsection}{3}\z@{.8\linespacing\@plus.7\linespacing}{.5\linespacing}{\normalfont\bfseries \centering}}
\makeatother

\numberwithin{equation}{section}

\newtheorem{theorem}{Theorem}[section]
\newtheorem{lemma}[theorem]{Lemma}
\newtheorem{proposition}[theorem]{Proposition}
\newtheorem{corollary}[theorem]{Corollary}
\newtheorem{introtheorem}{Theorem}

\theoremstyle{definition}
\newtheorem{definition}[theorem]{Definition}
\theoremstyle{remark}
\newtheorem{example}[theorem]{Example}
\newtheorem{remark}[theorem]{Remark}

\DeclareMathOperator{\Hom}{Hom}
\DeclareMathOperator{\Ext}{Ext}
\DeclareMathOperator{\rad}{rad}

\newcommand{\UU}{\mathcal U}
\newcommand{\Om}{\Omega}
\newcommand{\transpose}{^{\mathsf T}}
\newcommand{\kk}{\Bbbk}
\newcommand{\kdim}{\dim_{\kk}}
\newcommand{\Umod}[2]{U(#1,#2)}
\newcommand{\onevec}{\mathbf1}

\newcommand{\lemref}[1]{\hyperref[#1]{Lemma~\ref*{#1}}}
\newcommand{\propref}[1]{\hyperref[#1]{Proposition~\ref*{#1}}}
\newcommand{\thmref}[1]{\hyperref[#1]{Theorem~\ref*{#1}}}
\newcommand{\corref}[1]{\hyperref[#1]{Corollary~\ref*{#1}}}
\newcommand{\defref}[1]{\hyperref[#1]{Definition~\ref*{#1}}}
\newcommand{\exref}[1]{\hyperref[#1]{Example~\ref*{#1}}}
\newcommand{\figref}[1]{\hyperref[#1]{Figure~\ref*{#1}}}

\title[Syzygy matrices and higher extensions]{Syzygy matrices and higher extensions over string algebras}
\author{Martin Irungu}
\address{Ruhr-Universit\"at Bochum}
\email{martin.irungu@rub.de}
\date{}

\begin{document}
\begin{abstract}
Over a finite-dimensional string algebra, higher extension dimensions may eventually vanish, become periodic, or grow without bound. We show how their behaviour can be read from a finite directed graph. Second syzygies split into uniserial modules from a finite family; their syzygies define a nonnegative integer matrix, so higher extension dimensions become weighted counts of walks. Only the part of the graph relevant to the chosen pair of modules contributes. This gives four growth regimes: eventual vanishing, periodicity, quasi-polynomial growth, and exponential growth. We derive the construction from projective covers, give explicit extension formulas over truncated cycles, and illustrate the growth regimes with examples.
\end{abstract}
\maketitle

\tableofcontents

\section{Introduction}

For gentle algebras, higher extensions have a rather simple long-term behaviour. Baur and Schroll showed that, over an algebraically closed field, the extension dimensions for any two finite-dimensional modules eventually either vanish or repeat periodically \cite{BaurSchroll}. Their description also has a geometric interpretation on the surfaces associated with gentle algebras. Extensions of degree one had already been described rather explicitly, in terms of strings and their overlaps \cite{CPS}, and, for algebras from triangulations of unpunctured marked surfaces, in terms of crossings of curves \cite{CanakciSchroll17}. It is natural to ask what remains of this picture for the larger class of string algebras.

The answer is not simply more complicated periodicity. Polynomial and exponential growth can occur, sometimes over very small quivers. For example, a radical-square-zero algebra with two loops at one vertex has \(\dim_{\kk}\Ext_A^n(S,S)=2^n\), whereas two looped vertices joined by one arrow already give \(\dim_{\kk}\Ext_A^n(M,N)=n+1\). We will see both examples below. What is perhaps less obvious is that the behaviour depends on the target as well as on the source: a projective resolution may grow, while its extension spaces with a particular module stay bounded or even vanish.

One reason to study these numbers is that they measure how a projective resolution grows. For \(T=A/\rad A\), the number \(\kdim\Ext_A^n(M,T)\) is exactly the number of indecomposable projectives in its \(n\)-th term. A general target \(N\) sees only part of this resolution. In \exref{ex:polynomial}, for instance, a single source gives the three sequences \(1\), \(n\), and \(\binom n2\), depending on the target. In the bounded derived category these are the dimensions of the morphism spaces from \(M\) to the shifts \(N[n]\). Thus even when the modules themselves have a concrete description, their higher extensions can have quite different sizes.

Computing the first few terms may suggest a pattern, but it does not say when that pattern begins, whether it persists, or which part of the algebra is responsible for it. A finite graph answers these questions without requiring us to construct an infinite resolution. It also separates two sources of growth that will recur in the examples: a syzygy can split repeatedly within one family, or it can keep sending new summands into another.

For string algebras, the work of Huisgen-Zimmermann and Huisgen-Zimmermann--Smal\o\ provides the key observation: after two syzygies one meets only uniserial string modules, and there are finitely many possibilities \cite{HZ,HZS}. Our main results are \thmref{thm:intro-matrix} and \thmref{thm:intro-growth} below. The first gives a finite matrix description of higher extensions; the second classifies their eventual behaviour.

\begin{introtheorem}[the extension formula; \thmref{thm:matrix}]\label{thm:intro-matrix}
For every finite-dimensional string algebra \(A\) over a field \(\kk\), there is a finite nonnegative integer matrix \(B\), depending only on \(A\), with every row sum at most two. For each pair of finite-dimensional modules \(M,N\), there are nonnegative integer column vectors \(v_M,h_N\) such that
\begin{equation}\label{eq:intro}
 \kdim\Ext_A^i(M,N)=v_M\transpose B^{i-3}h_N\qquad(i\ge3).
\end{equation}
The vector \(v_M\) records the non-projective summands of \(\Om^2M\), and the entries of \(h_N\) are the first extension dimensions from the states into \(N\).
\end{introtheorem}

The matrix formula expresses extension dimensions as weighted counts of walks. Analysing these walks gives the second main result, the extension-dimension form of \thmref{thm:growth}.

\begin{introtheorem}[growth alternatives]\label{thm:intro-growth}
Let \(A\) be a finite-dimensional string algebra over a field \(\kk\), and let \(M,N\) be finite-dimensional modules. For \(a_i=\kdim\Ext_A^i(M,N)\), exactly one of the following holds:
\begin{enumerate}
\item \(a_i\) is eventually zero;
\item \(a_i\) is eventually periodic and not eventually zero;
\item \(a_i\) is eventually quasi-polynomial of positive degree;
\item \(\displaystyle\limsup_{i\to\infty}a_i^{1/i}=\rho\) for some \(1<\rho\le2\).
\end{enumerate}
\end{introtheorem}

Here quasi-polynomial means that, on each residue class modulo a fixed integer, the sufficiently large terms agree with a polynomial; the degree is the largest of their degrees. The graph explains which case occurs. When its cyclic components are directed cycles, staying in one gives periodic behaviour and passing through several can give polynomial growth. Branching within a strongly connected component gives exponential growth. Only the part of the graph which connects the two vectors matters.

There is a related connection with automata in the representation theory of string algebras. Rees constructs finite automata whose accepted words describe strings \cite{Rees08}. In the graph used here, a step instead records a syzygy summand. The extension formula is a weighted-automaton calculation, in the sense of \cite[\S1.6]{BerstelReutenauer11}.

The use of a syzygy quiver goes back to Cibils and was developed further by Howard, especially in the study of complexity and growth of projective resolutions \cite{Cibils93,Howard10,Howard15}. This article presents partial results from the author's already submitted dissertation \cite{IrunguDissertation}, with the truncated-cycle computations written in a more general form. The surface interpretation is treated there.

The article is organised as follows. Section~\ref{sec:prelim} recalls the homological background. Section~\ref{sec:matrix} sets up the syzygy matrix and proves the extension formula in \thmref{thm:matrix}. In Section~\ref{sec:growth}, we use its graph to distinguish exactly four growth patterns (\thmref{thm:growth}), followed by small examples in Section~\ref{sec:small-examples}. Section~\ref{sec:periodic} treats the periodic cases: gentle algebras and truncated cycles, where \thmref{thm:all-cyclic-ext} computes every extension dimension by counting eligible depths in the target. Finally, Section~\ref{sec:linked} shows how links between quiver cycles produce polynomial growth, with its degree determined by the accessible links.

\medskip
\noindent\textbf{Use of AI.} OpenAI's GPT-6 was used to improve the writing of this paper.

\section{Preliminaries}\label{sec:prelim}

Throughout, \(\kk\) is a field, \(Q\) is a finite quiver and \(A=\kk Q/I\), where \(J^m\subseteq I\subseteq J^2\) for some \(m\ge2\) and \(J\) is the arrow ideal. Paths are composed from left to right; all modules are finite-dimensional right modules. Write \(P_x=e_xA\), \(S_x=P_x/\rad P_x\), and \(T=A/\rad A=\bigoplus_{x\in Q_0}S_x\). In particular, each vertex simple is one-dimensional over \(\kk\); no assumption that \(\kk\) is algebraically closed is needed.

\begin{definition}[projective covers and syzygies]\label{def:covers}
A \emph{projective cover} of \(X\) is a surjection \(P(X)\twoheadrightarrow X\) from a projective module whose kernel lies in \(\rad P(X)\). Its kernel is the \emph{syzygy} \(\Om X\), giving the cover sequence
\[
0\longrightarrow\Om X\longrightarrow P(X)\longrightarrow X\longrightarrow0.
\]
Iterating gives the \emph{minimal projective resolution} of \(X\).
\end{definition}

Every module has a projective cover, unique up to isomorphism, and a decomposition into indecomposables, unique up to order and isomorphism \cite[I, Theorem~4.2; II, Theorem~2.2]{ARS}. We recall the standard consequences for syzygies and extensions below; see \cite[I, Proposition~4.3 and \S5]{ARS}.

\begin{lemma}[Projective covers and extensions]\label{lem:homological}
For finite-dimensional \(A\)-modules \(X,Y\) the following hold.
\begin{enumerate}
\item Syzygies are additive: \(\Om(X\oplus Y)\cong\Om X\oplus\Om Y\); projective summands disappear. In particular, \(\Ext_A^i(P,Y)=0\) for \(i\ge1\) if \(P\) is projective.
\item For \(i\ge1\), dimension shift gives
\[
\Ext_A^{i+1}(X,Y)\cong\Ext_A^i(\Om X,Y).
\]
Moreover,
\begin{equation}\label{eq:cokernel}
\Ext_A^1(X,Y)\cong\operatorname{coker}\bigl(\Hom_A(P(X),Y)\longrightarrow\Hom_A(\Om X,Y)\bigr).
\end{equation}
Consequently,
\[
\kdim\Ext_A^1(X,Y)=\kdim\Hom_A(\Om X,Y)-\kdim\Hom_A(P(X),Y)+\kdim\Hom_A(X,Y).
\]
\item If \(\cdots\to P_1\to P_0\to X\) is a minimal projective resolution, then
\[
\Ext_A^i(X,T)\cong\Hom_A(P_i,T)\qquad(i\ge0).
\]
Thus \(\kdim\Ext_A^i(X,T)\) counts the indecomposable summands of \(P_i\), with multiplicity.
\end{enumerate}
\end{lemma}

For string modules we follow the construction in \cite[\S3]{ButlerRingel}. The classification over an arbitrary field is also given in \cite[Theorems~1.1--1.2]{CrawleyBoevey18}. We recall the part of the construction needed for the cover calculations.

\section{The syzygy matrix}\label{sec:matrix}

\subsection{Strings and projective arms}
Write \(s(\alpha)\) and \(t(\alpha)\) for the source and target of an arrow \(\alpha\).

\begin{definition}[string and gentle algebras]\label{def:algebras}
The algebra \(A=\kk Q/I\) is a \emph{string algebra} if \(I\) is generated by paths, at most two arrows enter and at most two leave each vertex, and each arrow has at most one nonzero successor and one nonzero predecessor. More explicitly, for each \(\alpha\) there is at most one composable \(\beta\) with \(\alpha\beta\notin I\) and at most one composable \(\gamma\) with \(\gamma\alpha\notin I\).

It is \emph{gentle} if, in addition, \(I\) is generated by paths of length two and each arrow has at most one composable zero successor and one composable zero predecessor.
\end{definition}

Throughout this section, \(A\) is a string algebra.

\begin{definition}[strings and string modules]\label{def:strings}
A \emph{basis path} is a path with nonzero image in \(A\). A \emph{string} is a composable word in arrows and formal inverses, without immediate cancellation or a subword belonging to \(I\) or its inverse. Trivial strings are allowed; a string containing only arrows is \emph{direct}.

The module \(M(w)\) has one basis vector at each successive occurrence of a vertex in \(w\); an arrow in the word acts between the corresponding vectors in its own direction. All other arrow actions are zero. A \emph{peak} is a basis occurrence that survives in the top; a \emph{valley} is one in the socle.
\end{definition}

Occurrences remain distinct even when their vertex labels agree, and \(M(w)\cong M(w^{-1})\). The projective cover of a string module has one summand for each peak, so a string with one peak is cyclic. A direct string \(p\) gives a uniserial module \(M(p)\), whose submodules are totally ordered. Every uniserial string module has a direct representative after reversing the word if necessary.

\begin{definition}[projective arms]\label{def:arms}
For each arrow \(\beta\) from \(x\), the unique maximal basis path beginning with \(\beta\) is the \(\beta\)-\emph{arm} of \(P_x\).
\end{definition}

The path basis of \(P_x\) consists of \(e_x\) and the nontrivial prefixes of its at most two arms. Distinct arms have disjoint path bases, even if they visit the same quiver vertices.

\begin{lemma}[the uniserial kernel]\label{lem:arms}
Let \(p\) be a direct basis path starting at \(x\). The cover \(P_x\twoheadrightarrow M(p)\) retains the initial segment specified by \(p\) on its arm and kills every other arm. Its kernel is the direct sum of the remaining tails. In particular, \(\Om M(p)\) is a sum of at most two uniserial string modules.
\end{lemma}
\begin{proof}
The cover sends \(e_x\) to the top generator. The prefixes of \(p\) map to the successive basis vectors of \(M(p)\); all other basis paths are killed. On each arm, the killed paths form a final segment, generated by its first path. Unique continuation makes this segment uniserial. The segments have disjoint path bases and together span the kernel. For a trivial \(p\), both arms are killed and the same argument applies.
\end{proof}

\begin{figure}[ht]
\centering
\begin{tikzpicture}[scale=.95]
\node (x) at (-2.8,0) {\(e_x\)};
\coordinate (a) at (-1.8,.7); \coordinate (b) at (-.8,.7); \coordinate (c) at (.2,.7);
\coordinate (d) at (1.2,.7); \coordinate (e) at (2.2,.7);
\coordinate (f) at (-1.8,-.7); \coordinate (g) at (-.8,-.7);
\fill (a) circle (2pt); \fill (b) circle (2pt); \fill (c) circle (2pt);
\draw[fill=white] (d) circle (2.5pt); \draw[fill=white] (e) circle (2.5pt);
\draw[fill=white] (f) circle (2.5pt); \draw[fill=white] (g) circle (2.5pt);
\draw[->] (x)--(a) node[midway,above] {\(\beta_1\)};
\draw[->] (a)--(b); \draw[->] (b)--(c);
\draw[->,dashed] (c)--(d); \draw[->] (d)--(e);
\draw[->,dashed] (x)--(f) node[midway,below left] {\(\beta\)}; \draw[->] (f)--(g);
\node at (-.75,1.05) {\(M(p)\)}; \node at (1.65,1.05) {\(T_1\)}; \node at (-1.3,-1.05) {\(T_2\)};
\end{tikzpicture}
\caption{A uniserial cover and its kernel tails.}
\label{fig:arms}
\end{figure}
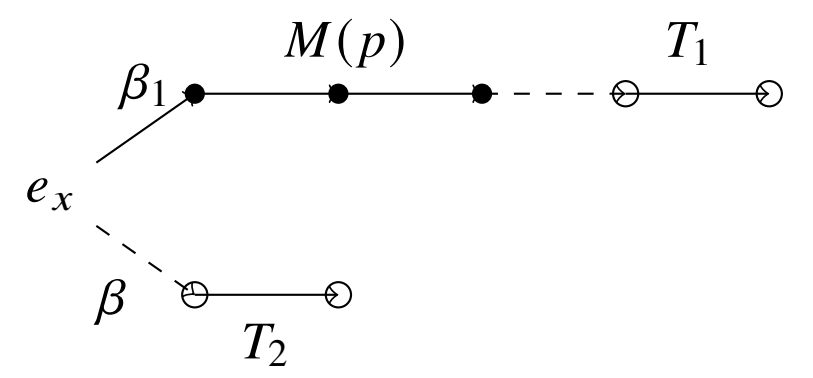

In \figref{fig:arms}, filled dots survive in the quotient and open dots form the kernel. An exhausted arm contributes zero, not a trivial string; a one-vector tail contributes a simple. The kernel vanishes exactly when \(p\) exhausts the sole arm (or \(p\) is trivial at a sink).

\begin{example}[why arrow ideals alone need not suffice]\label{ex:long-tail}
Consider the bound quiver
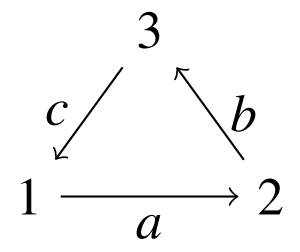
\[
\begin{tikzpicture}[baseline=(current bounding box.center)]
\node (a0) at (0,0) {\(1\)}; \node (a1) at (1.6,0) {\(2\)}; \node (a2) at (.8,1.1) {\(3\)};
\draw[->] (a0)--(a1) node[midway,below] {\(a\)};
\draw[->] (a1)--(a2) node[midway,right] {\(b\)};
\draw[->] (a2)--(a0) node[midway,left] {\(c\)};
\end{tikzpicture}
\qquad I=J^3.
\]
The ideal \(aA\cong M(b)\) has cover \(P_{t(a)}\) with path basis \(e_{t(a)},b,bc\). Its kernel is the simple \(\kk\,bc\). Every arrow ideal has dimension two, so this kernel is not an arrow ideal. The uniserial state set must also contain simples.
\end{example}


\subsection{Reduction after two syzygies}
Set \(\Om^0M=M\), and write \(\underline X\) for the non-projective part of a module \(X\). The reason a finite matrix suffices is a theorem of Huisgen-Zimmermann for monomial algebras. We use its following consequence for string algebras.

\begin{proposition}[second-syzygy reduction]\label{prop:reduction}
Every second syzygy is a direct sum of ideals \(pA\), where \(p\) is a positive-length basis path. In particular, it is a direct sum of uniserial string modules from a fixed finite family.
\end{proposition}
\begin{proof}
The monomial second-syzygy theorem \cite[Theorem~I]{HZ}, in the form \cite[Theorem~A(1)]{HZDomino}, says that the syzygy of a submodule of a projective is a direct sum of principal ideals generated by positive-length paths. Apply it to \(\Om M\subseteq P(M)\), passing to the opposite algebra to obtain the right-module statement. The theorem holds over an arbitrary field.

For a string algebra, a positive-length path has at most one nonzero continuation at each step. The path basis of \(pA\) therefore forms a single chain, so \(pA\) is a uniserial string module. There are only finitely many basis paths.
\end{proof}

For a string module, we can see the first kernel directly in its cover. The description is that of \cite[Proposition~1(1)(a)]{HZS}; the path-basis calculation below also works over an arbitrary field.

There is one projective for each peak of the string. Its basis paths follow the string until they reach an end or a valley; beyond that point they map to zero. At an interior valley, two paths \(p,q\) reach the same basis vector, so their difference \(\zeta=(p,-q)\) belongs to the kernel. Their final arrows are distinct. An outgoing arrow can therefore extend at most one of them, and every nonzero extension of either path belongs to \(\zeta A\).

The kernel is consequently the direct sum of the unused endpoint tails and these valley modules. For a trivial string, the two possible arms are its endpoint tails. They exhaust the killed paths and the differences, and have disjoint supports in the path basis of the cover. Here, as in the string itself, repeated vertex labels still represent distinct basis occurrences. Each summand has one peak, so its own cover has only arm tails in its kernel. This gives the second-syzygy reduction for strings. A valley module need not itself be projective.

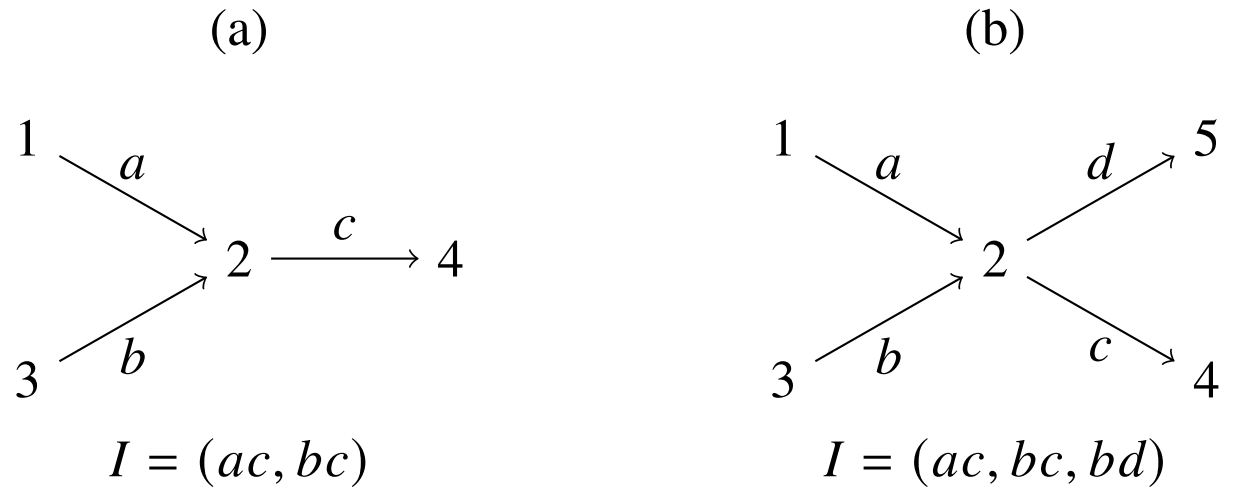
\begin{figure}[ht]
\centering
\begin{tikzpicture}
\begin{scope}
\node at (1.4,1.5) {(a)};
\node (v1) at (0,.8) {\(1\)}; \node (v3) at (0,-.8) {\(3\)};
\node (v2) at (1.4,0) {\(2\)}; \node (v4) at (2.8,0) {\(4\)};
\draw[->] (v1)--(v2) node[midway,above] {\(a\)};
\draw[->] (v3)--(v2) node[midway,below] {\(b\)};
\draw[->] (v2)--(v4) node[midway,above] {\(c\)};
\node at (1.4,-1.35) {\(I=(ac,bc)\)};
\end{scope}
\begin{scope}[shift={(5,0)}]
\node at (1.4,1.5) {(b)};
\node (v1) at (0,.8) {\(1\)}; \node (v3) at (0,-.8) {\(3\)};
\node (v2) at (1.4,0) {\(2\)}; \node (v4) at (2.8,-.8) {\(4\)}; \node (v5) at (2.8,.8) {\(5\)};
\draw[->] (v1)--(v2) node[midway,above] {\(a\)};
\draw[->] (v3)--(v2) node[midway,below] {\(b\)};
\draw[->] (v2)--(v4) node[midway,below] {\(c\)};
\draw[->] (v2)--(v5) node[midway,above] {\(d\)};
\node at (1.4,-1.35) {\(I=(ac,bc,bd)\)};
\end{scope}
\end{tikzpicture}
\caption{Valley kernels: (a) \(S_2\); (b) \(M(d)\).}
\label{fig:valleys}
\end{figure}

\begin{example}[a valley kernel need not be projective]\label{ex:valley}
Consider the two bound quivers in \figref{fig:valleys}. In both, the cover of \(M(ab^{-1})\) is \(P_1\oplus P_3\), and the two paths into the valley give the kernel generator \(\zeta=(a,-b)\). In (a), all arrows kill \(\zeta\), so
\[
\Om M(ab^{-1})\cong S_2,\qquad \dim P_2=2.
\]
In (b), \(\zeta d=(ad,0)\ne0\), and the kernel is \(M(d)\). The extra killed path \(ad\) is already part of the valley module. Both algebras are string but not gentle: \(c\) has two zero predecessors.
\end{example}

\subsection{The matrix and the two endpoint vectors}
Let \(\UU\) contain one representative of each isomorphism class of non-projective uniserial \emph{string} modules. It is finite because there are only finitely many basis paths. \lemref{lem:arms} shows that taking syzygies and deleting projective summands stays within direct sums of this family.

For an indecomposable \(V\), write \((X:V)\) for its multiplicity in \(X\).

\begin{definition}[the syzygy matrix and endpoint vectors]\label{def:syzygy-matrix}
The \emph{syzygy matrix} \(B\) on \(\UU\) and the vectors associated with \(M,N\) are
\[
\begin{aligned}
 B_{UV}&=(\underline{\Om U}:V),\\
 v_M&=\big((\underline{\Om^2M}:U)\big)_{U\in\UU},&
 h_N&=\big(\kdim\Ext_A^1(U,N)\big)_{U\in\UU}.
\end{aligned}
\]
Rows index source modules; columns index their syzygy summands. We regard \(v_M,h_N\) as column vectors and write \(v\transpose\) for the transpose of \(v\).
\end{definition}

The matrix depends on \(A\), while the two vectors record the chosen pair.

%

\begin{theorem}[the extension formula]\label{thm:matrix}
Every row of \(B\) has sum at most two. For all \(M,N\),
\[
 \kdim\Ext_A^i(M,N)=v_M\transpose B^{i-3}h_N\qquad(i\ge3).
\]
\end{theorem}
\begin{proof}
By \lemref{lem:arms}, each state has at most two syzygy summands, giving the row bound. By \propref{prop:reduction}, all non-projective summands of \(\Om^2M\) belong to \(\UU\).

Let \(m_n\) be the column with entries \(m_n(U)=(\underline{\Om^{n+2}M}:U)\), so \(m_0=v_M\). Each copy of \(U\) contributes \(B_{UV}\) copies of \(V\) at the next step. Additivity of syzygy, and the fact that projective summands disappear, therefore give
\[
 m_{n+1}(V)=\sum_{U\in\UU}m_n(U)B_{UV},\qquad m_{n+1}\transpose=m_n\transpose B.
\]
Induction yields \(m_n\transpose=v_M\transpose B^n\). Now shift an extension of degree \(i\ge3\) down to degree one:
\[
\begin{aligned}
 \kdim\Ext_A^i(M,N)&=\kdim\Ext_A^1(\Om^{i-1}M,N)\\
 &=\sum_{U\in\UU}m_{i-3}(U)\,\kdim\Ext_A^1(U,N)\\
 &=v_M\transpose B^{i-3}h_N.
\end{aligned}
\]
Projective summands of \(\Om^{i-1}M\) contribute zero to this sum.
\end{proof}

The exponent begins at zero in degree three. The first two positive degrees require their own cover calculations. For \(h_N\), use the cokernel or the dimension formula in \lemref{lem:homological}; maps \(\Om U\to N\) that extend to the cover do not contribute.

\begin{example}[a complete calculation]\label{ex:truncated-loop}
Let \(A=\kk[x]/(x^3)\), \(S=A/(x)\) and \(U=A/(x^2)\). The two covers give \(\Om S\cong U\) and \(\Om U\cong S\). In the order \((S,U)\),
\[
 B=\begin{pmatrix}0&1\\1&0\end{pmatrix}.
\]
Take \(M=S\oplus U^{\oplus2}\oplus A\) and \(N=S\). Then \(\Om^2M\cong S\oplus U^{\oplus2}\), so \(v_M=(1,2)\transpose\). Every map \(A\to S\) kills \(\rad A\); both restriction maps in \eqref{eq:cokernel} are zero. Since \(\Hom_A(S,S)\) and \(\Hom_A(U,S)\) each have dimension one, \(h_S=(1,1)\transpose\). Hence
\[
 \kdim\Ext_A^i(M,S)=(1,2)B^{i-3}(1,1)\transpose=3\qquad(i\ge3).
\]
The same covers give \(3\) in degrees one and two. Degree zero is \(4\): the projective summand contributes to Hom, but not to positive Ext.
\end{example}

Once \(B\) is known, a new pair of modules requires only the two vectors \(v_M\) and \(h_N\). For gentle string pairs over an algebraically closed field, the arrow and overlap basis of \cite[Theorem~A]{CPS} is another way to compute \(h_N\).

\begin{corollary}[rational series]\label{cor:series}
The truncated Poincar\'e series, regarded as a formal power series over \(\mathbb Q\), is
\[
 \sum_{i\ge3}\kdim\Ext_A^i(M,N)t^i
   =t^3v_M\transpose(I-tB)^{-1}h_N.
\]
It is rational, with denominator dividing \(\det(I-tB)\). The higher extension dimensions satisfy a linear recurrence.
\end{corollary}
\begin{proof}
Sum the matrix powers as a formal geometric series and use the adjugate formula. A recurrence follows by applying the characteristic polynomial of \(B\) to its powers.
\end{proof}


The denominator can be smaller for a particular pair. In \exref{ex:truncated-loop}, \(\det(I-tB)=1-t^2\), but the truncated series is \(3t^3/(1-t)\). The syzygies alternate while the target gives both states the same weight. The constant tail \(3\) needs only a one-state matrix description, \(3\cdot(1)^n\cdot1\). In general, linear minimization over \(\mathbb Q\) can reduce the number of states further \cite[\S2.3]{BerstelReutenauer11}, though the resulting states need not represent modules or have nonnegative integral transition weights.

\section{Growth from the syzygy graph}\label{sec:growth}

View a nonnegative integer matrix \(B\) as the adjacency matrix of a directed multigraph, with \(B_{UV}\) edges from \(U\) to \(V\). Thus \((B^n)_{UV}\) counts walks of length \(n\), and
\[
 a_n=v\transpose B^nh
\]
counts them with initial and terminal weights. This is a weighted automaton over \(\mathbb N\) with one letter: each reading of that letter applies \(B\), and the word of length \(n\) has weight \(a_n\). We use the same walk-counting methods that appear in the growth theory of regular languages \cite{GKRS10}.

\begin{definition}[the trimmed graph]\label{def:trimmed}
The \emph{support} of a vector \(v\) is \(\operatorname{supp}v=\{U:v_U\ne0\}\). A vertex is \emph{accessible} if it can be reached from \(\operatorname{supp}v\), and \emph{coaccessible} if it can reach \(\operatorname{supp}h\); walks of length zero are allowed. Keeping only vertices with both properties gives the \emph{trimmed graph} for \((v,h)\), or the \emph{trimmed syzygy graph} for \((M,N)\) when \(v=v_M\) and \(h=h_N\).
\end{definition}

This is the terminology of automata theory \cite{Sakarovitch}. Restricting to the trimmed graph does not change \(a_n\).

\begin{definition}[components and periods]\label{def:components}
A \emph{strongly connected component} is a maximal set of vertices each reachable from every other. It is \emph{cyclic} if it contains a cycle, and it is a \emph{directed cycle} when its internal edges, counted with multiplicity, form a single cycle (possibly a loop). The \emph{period} of a cyclic component is the greatest common divisor of the lengths of its positive-length closed walks.
\end{definition}

Collapsing the components to single vertices gives an acyclic graph. Ordering them so that edges point forward puts the adjacency matrix in block triangular form.

\begin{definition}[spectral terminology]\label{def:spectral}
The \emph{spectral radius} \(\rho(C)\) of a square matrix \(C\) is the largest absolute value of its eigenvalues. A nonnegative matrix is \emph{irreducible} if its graph is strongly connected, and \emph{primitive} if some positive power has every entry positive. When \(\rho(C)\) is itself an eigenvalue, it is called the \emph{Perron eigenvalue}. A positive \emph{left Perron eigenvector} is a column \(y>0\) satisfying \(y\transpose C=\rho(C)y\transpose\); a right one satisfies \(Cx=\rho(C)x\).
\end{definition}

For an irreducible nonnegative matrix \(C\ne[0]\), Perron--Frobenius theory gives \(\rho(C)>0\) and positive right and left Perron eigenvectors. We use these facts and the decomposition into periodic classes from \cite[Chapter~1]{SenetaNonnegative}. Write \(\onevec\) for the column with every entry one.

\subsection{The four alternatives}
\begin{definition}[eventual growth]\label{def:growth}
A sequence \((a_n)\) is \emph{eventually zero} if all sufficiently large terms vanish, and \emph{eventually periodic} if \(a_{n+L}=a_n\) for some \(L>0\) and all sufficiently large \(n\). It is \emph{eventually quasi-polynomial of degree \(d\)} if, on each residue class modulo some \(L>0\), its sufficiently large terms agree with a polynomial of degree at most \(d\), with degree exactly \(d\) on at least one class. Such an \(L\) is a \emph{quasiperiod}.
\end{definition}

Degree zero means eventual periodicity with a nonzero tail.

\begin{lemma}\label{lem:radius-one}
Let \(C\) be the adjacency matrix of a cyclic strongly connected directed multigraph. Then \(\rho(C)\ge1\), with equality precisely when the graph is a directed cycle.
\end{lemma}
\begin{proof}
Every row sum is at least one. If \(y>0\) is a left Perron eigenvector, then \(y\transpose C\onevec\ge y\transpose\onevec\), with strict inequality whenever a row sum exceeds one. Thus \(\rho(C)\ge1\), with equality only when every row sum is one. Strong connectivity then gives a directed cycle \cite[Chapter~1]{SenetaNonnegative}.
\end{proof}

\begin{theorem}[growth alternatives]\label{thm:growth}
Let \(a_n=v\transpose B^nh\) for a nonnegative integer matrix \(B\) and nonnegative integer vectors \(v,h\). Exactly one of the following holds.
\begin{enumerate}
\item The trimmed graph is acyclic. Then \(a_n\) is eventually zero.
\item The trimmed graph contains a cycle, every cyclic component is a directed cycle, and no walk meets two cyclic components. Then \(a_n\) is eventually periodic and not eventually zero.
\item Every cyclic component is a directed cycle, and the greatest number met by a walk is \(r\ge2\). Then \(a_n\) is eventually quasi-polynomial of degree \(r-1\).
\item Some cyclic component of the trimmed graph has spectral radius greater than one. If \(\rho\) is the largest such radius, then
\[
\limsup_{n\to\infty}a_n^{1/n}=\rho>1.
\]
\end{enumerate}
In \textup{(2)} and \textup{(3)}, the least common multiple of the cycle lengths is a quasiperiod. The same conclusions hold for higher extension dimensions, with the index shift in \thmref{thm:matrix}.
\end{theorem}
\begin{proof}
An acyclic graph has bounded walk length. Every cycle in the trimmed graph lies on a walk from the initial to the terminal support, so it gives arbitrarily long counted walks. Suppose now that all cyclic components are directed cycles. Fix the route through the component graph, its connecting edges, and its initial and terminal support vertices. In each cycle, a walk follows the unique path from entry to exit of length less than one circuit, with any number of complete circuits inserted. There are finitely many such routes \(\gamma\), and each contributes
\[
 \frac{c_\gamma t^{d_\gamma}}{\prod_{j=1}^{s_\gamma}(1-t^{\ell_{\gamma,j}})},
 \qquad c_\gamma>0,
\]
to \(\sum_{n\ge0}a_nt^n\), where \(s_\gamma\) is the number of cycles on the route. Routes without cycles contribute a polynomial.

Let \(L\) be the least common multiple of the cycle lengths. All poles of the series are \(L\)-th roots of unity, so its coefficients are eventually quasi-polynomial with quasiperiod \(L\) \cite[Proposition~4.4.1]{StanleyEC1}. To find the degree, put \(r=\max_\gamma s_\gamma\). The pole at \(t=1\) has order exactly \(r\), since the positive route contributions cannot cancel there, and no other pole has larger order. In a partial fraction expansion, a pole of order \(j\) contributes a polynomial of degree \(j-1\) times an \(L\)-periodic sequence. The degree is therefore \(r-1\). The cases \(r=1\) and \(r\ge2\) give \textup{(2)} and \textup{(3)}.

Finally, let \(D\) be the adjacency matrix of a cyclic component of maximal radius \(\rho>1\). The trimmed matrix is block triangular with the component matrices on its diagonal, so its spectral radius is \(\rho\). The spectral-radius formula for matrix powers gives \(\limsup a_n^{1/n}\le\rho\).

For the reverse inequality, fix entry and exit walks through a vertex \(z\) of \(D\). They give \(a_{n+c}\ge K(D^n)_{zz}\) for some \(c\) and \(K>0\). If \(d\) is the period of \(D\), the block of \(D^d\) containing \(z\) is primitive and has spectral radius \(\rho^d\). Perron--Frobenius theory therefore gives \((D^{dm})_{zz}\sim c_z\rho^{dm}\), with \(c_z>0\). This subsequence proves the reverse inequality and hence \textup{(4)}; \lemref{lem:radius-one} exhausts the possibilities.
\end{proof}

The row bound gives \(\kdim\Ext_A^i(M,N)\le \lVert v_M\rVert_1\lVert h_N\rVert_\infty 2^{i-3}\) for \(i\ge3\), so an exponential rate is at most two. Trimming matters: the full matrix may grow faster than the coefficient for a chosen pair. Even in the exponential case some residue classes can vanish.



\begin{remark}[finer asymptotics]
A limit need not exist: residue classes may vanish or have different rates. Finer asymptotics follow from periodic-class decomposition and Perron--Frobenius theory \cite{Rothblum75,SenetaNonnegative}. For example, \(B=\left(\begin{smallmatrix}0&2\\1&0\end{smallmatrix}\right)\) and \(v=h=(1,0)\transpose\) give \(a_{2m}=2^m\) and \(a_{2m+1}=0\). Thus \(\limsup a_n^{1/n}=\sqrt2\), although every odd term vanishes.
\end{remark}



\subsection{Small examples of the four behaviours}\label{sec:small-examples}
For the following examples we impose all paths of length two as relations. Their extensions can be computed directly from the quiver.

\begin{lemma}[Radical-square-zero simples]\label{lem:radical2}
For \(A=\kk Q/J^2\), let \(D_{xy}\) be the number of arrows from \(x\) to \(y\). Then
\[
\Om S_x\cong\bigoplus_y S_y^{\oplus D_{xy}},\qquad
\kdim\Ext_A^i(S_x,S_y)=(D^i)_{xy}\quad(i\ge0).
\]
\end{lemma}
\begin{proof}
The radical of \(P_x\) is the direct sum of the simples at the targets of arrows from \(x\). Iterate and apply \lemref{lem:homological}.
\end{proof}

\begin{figure}[ht]
\centering
\begin{tikzpicture}[scale=.93]
\begin{scope}[shift={(0,2.0)}]
\node at (1.5,1) {(a)};
\node (a) at (0,.2) {\(1\)}; \node (b) at (1.5,.2) {\(2\)}; \node (c) at (3,.2) {\(3\)};
\draw[->] (a)--(b); \draw[->] (b)--(c);
\end{scope}
\begin{scope}[shift={(5.0,2.0)}]
\node at (1.2,1) {(b)};
\node (d) at (0,.2) {\(1\)}; \node (e) at (2.4,.2) {\(2\)};
\draw[->] (d) to[bend left=20] (e); \draw[->] (e) to[bend left=20] (d);
\end{scope}
\begin{scope}
\node at (1.5,1) {(c)};
\node (f) at (0,.2) {\(0\)}; \node (g) at (1.5,.2) {\(1\)}; \node (h) at (3,.2) {\(2\)};
\draw[->] (f)--(g); \draw[->] (g)--(h);
\draw[->] (f) edge[loop above] (f); \draw[->] (g) edge[loop above] (g); \draw[->] (h) edge[loop above] (h);
\end{scope}
\begin{scope}[shift={(5.0,0)}]
\node at (1.2,1) {(d)};
\node (i) at (0,.2) {\(1\)}; \node (j) at (2.4,.2) {\(2\)};
\draw[->] (i) edge[loop above] (i);
\draw[->] (i) to[bend left=20] (j); \draw[->] (j) to[bend left=20] (i);
\end{scope}
\end{tikzpicture}
\caption{Four growth behaviours for \(\kk Q/J^2\).}
\label{fig:small}
\end{figure}
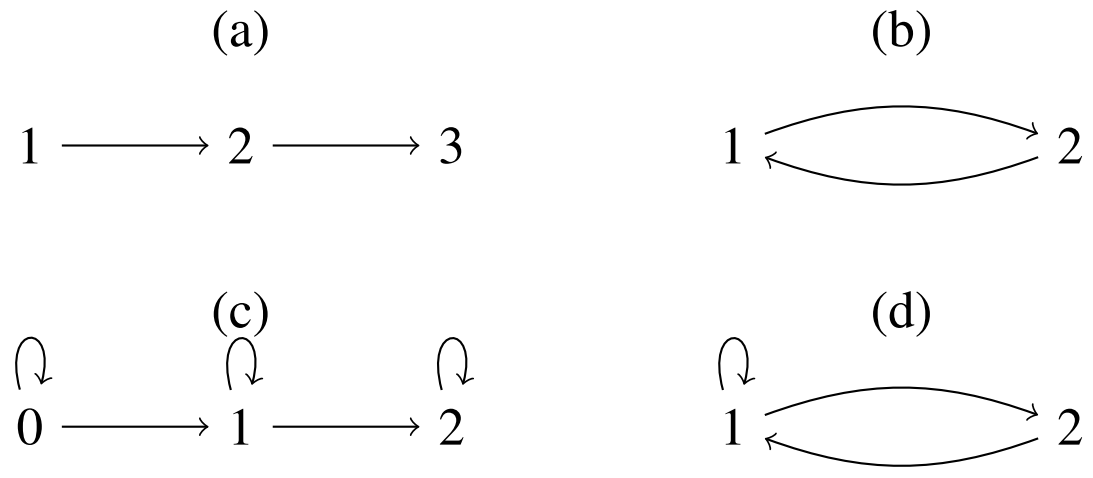

In \figref{fig:small}, each vertex has at most two incoming and two outgoing arrows, so these radical-square-zero algebras are string. A walk records successive syzygies; it need not survive as a nonzero path in \(A\).

\begin{example}[vanishing and periodicity]
In \figref{fig:small}(a), the only nonzero extension between \(S_1\) and \(S_3\) occurs in degree two, where it has dimension one. Although the length-two path is zero in \(A\), it still counts as a walk. In (b),
\[
  \kdim\Ext_A^n(S_1,S_1)=\begin{cases}1& n\text{ even},\\0&n\text{ odd}.\end{cases}
\]
\end{example}

\begin{example}[one or two transfers]\label{ex:polynomial}
In \figref{fig:small}(c), a walk from \(0\) to \(1\) chooses one of its \(n\) steps at which to cross the connecting arrow. A walk from \(0\) to \(2\) chooses two crossing steps. Thus
\[
 \begin{aligned}
 \kdim\Ext_A^n(S_0,S_0)&=1,\\
 \kdim\Ext_A^n(S_0,S_1)&=n,&
 \kdim\Ext_A^n(S_0,S_2)&=\binom n2.
 \end{aligned}
\]
The last two series are \(t/(1-t)^2\) and \(t^2/(1-t)^3\), respectively. The source is the same in all three cases; changing the target changes the growth.
\end{example}

\begin{example}[Fibonacci growth]\label{ex:fibonacci}
In \figref{fig:small}(d), the adjacency matrix is
\[
 F=\begin{pmatrix}1&1\\1&0\end{pmatrix}.
\]
The walk counts satisfy \(a_0=a_1=1\) and \(a_n=a_{n-1}+a_{n-2}\) for \(n\ge2\). Thus \(\kdim\Ext_A^n(S_1,S_1)\) is \(1,1,2,3,5,8,\ldots\), with generating function \(1/(1-t-t^2)\) and exponential rate \((1+\sqrt5)/2\). Here branching occurs within one recurrent component.
\end{example}

\begin{example}[two loops]\label{ex:two-loops}
For the bound quiver
\[
\begin{tikzpicture}[baseline=(v.base)]
\node (v) at (0,0) {\(0\)};
\draw[->] (v) edge[loop left] node[left] {\(x\)} (v);
\draw[->] (v) edge[loop right] node[right] {\(y\)} (v);
\end{tikzpicture}
\qquad I=(x^2,xy,yx,y^2),
\]
each step of a walk has two choices. Thus \(A=\kk\langle x,y\rangle/(x,y)^2\) has \(\kdim\Ext_A^n(S,S)=2^n\). The full non-projective uniserial string family is \((S,M(x),M(y))\), with matrix
\[
 B=\begin{pmatrix}2&0&0\\1&0&0\\1&0&0\end{pmatrix}.
\]
The projective \(A\) has radical \(S^{\oplus2}\); the covers of \(M(x)\) and \(M(y)\) each have simple kernel. Only \(S\) is recurrent; the two other states feed into it after one step. Quadratic relations alone therefore allow exponential growth.
\end{example}

\section{Periodic examples}\label{sec:periodic}

\subsection{The gentle case}
Write \(R(\alpha)=\alpha A\) for an arrow ideal. In the small gentle quiver below, the dotted curve marks the relation \(ab=0\):
\[
\begin{tikzpicture}[baseline=(current bounding box.center)]
\node (a) at (0,0) {\(1\)}; \node (b) at (2,0) {\(2\)}; \node (c) at (4,0) {\(3\)};
\draw[->] (a)--(b) node[midway,below] {\(a\)};
\draw[->] (b)--(c) node[midway,below] {\(b\)};
\draw[densely dotted] (1.1,.18) to[bend left=35] (2.9,.18);
\end{tikzpicture}
\qquad \Om R(a)\cong R(b).
\]

The gentle conditions simplify both parts of the kernel calculation. Quadraticity also makes every positive-length path ideal isomorphic to an arrow ideal. Together with \propref{prop:reduction}, this gives periodicity over an arbitrary field.

\begin{lemma}[gentle kernels]\label{lem:gentle}
For a gentle algebra, \(\Om R(\alpha)\) is \(R(\beta)\) if \(\alpha\beta\in I\) for an arrow \(\beta\) starting at \(t(\alpha)\), and is zero if no such arrow exists. For a string module,
\[
 \Om M(w)\cong L_0\oplus P\oplus L_1,
\]
where \(P\) is projective and each \(L_j\) is zero or an arrow ideal.
\end{lemma}
\begin{proof}
The cover \(P_{t(\alpha)}\to\alpha A\) sends \(u\) to \(\alpha u\). A basis path is killed precisely when it begins with an arrow \(\beta\) such that \(\alpha\beta\in I\): quadraticity excludes a first obstruction farther along the path. There is at most one such \(\beta\), so the kernel is \(\beta A\) or zero.

For the second assertion, use the string-kernel description following \propref{prop:reduction}. At a valley \(x\), the final arrows of \(p,q\) are distinct. Each outgoing arrow has exactly one nonzero predecessor among them, by the two gentle predecessor conditions. Quadraticity implies that every positive-length basis path from \(x\) extends exactly one of \(p,q\). Distinct paths have disjoint images in the path basis of the cover, so \(P_x\to(p-q)A\), \(e_x\mapsto p-q\), is injective as well as surjective.

For any positive-length basis path \(u=v\alpha\), quadraticity gives \(uw\ne0\) if and only if \(\alpha w\ne0\) for every basis path \(w\) starting at \(t(\alpha)\). Thus \(\alpha w\mapsto uw\) defines an isomorphism \(\alpha A\cong uA\). In particular, each endpoint tail is an arrow ideal up to isomorphism. This gives the decomposition of \cite[Proposition~3.1]{BaurSchroll}.
\end{proof}

\begin{corollary}[gentle periodicity]\label{cor:gentle}
Higher extension dimensions over a gentle algebra over any field are eventually periodic, including the possibility of eventual zero.
\end{corollary}
\begin{proof}
By \propref{prop:reduction} and the path-ideal isomorphism just proved, \(\Om^2M\) is a direct sum of arrow ideals. The first part of \lemref{lem:gentle} sends each arrow ideal to another arrow ideal or to zero. There are finitely many arrows, so each such sequence terminates or becomes periodic. Taking a common period for the finitely many summands makes the higher syzygies of \(M\) periodic. Dimension shift gives eventual periodicity of \(\kdim\Ext_A^i(M,N)\) for every target \(N\).
\end{proof}

Over an algebraically closed field this recovers the vanishing-or-periodicity conclusion of \cite[Theorem~3.4]{BaurSchroll}, which also gives a finer endpoint criterion for a nonzero periodic tail. In that setting, bands satisfy a stronger conclusion: their projective and injective dimensions are at most one. This is \cite[Corollary~3.6]{BaurSchroll} for quasi-simple bands; a filtration by quasi-simple bands and the long exact sequences for Ext give the same bounds in higher rank. For the radical-square-zero cycle \(\kk C_m/J^2\) one sees the period directly: \(\Om S_j=S_{j+1}\) and
\[
 \kdim\Ext_A^i(S_j,S_l)=[l\equiv j+i\pmod m].
\]
Brackets denote the indicator of the enclosed condition.

%

\subsection{Truncated cycles}
Let
\[
n\ge1,\qquad r\ge2,\qquad \Lambda_{n,r}=\kk C_n/J^r,
\]
where \(C_n\) has arrows \(a_j:j\to j+1\), with indices modulo \(n\). Let \(\Umod{j}{\ell}\) be the uniserial module with top \(S_j\) and length \(\ell\), for \(1\le\ell\le r\). Its successive composition factors are \(S_j,\ldots,S_{j+\ell-1}\), and \(U(j,r)=P_j\). We allow \(r>n\), so vertex labels can repeat in a composition series. The cover calculation is unchanged; the Hom calculation must count each eligible occurrence separately.

\begin{proposition}[complementary tails]\label{prop:cyclic}
For \(1\le\ell<r\),
\[
 \Om \Umod{j}{\ell}\cong \Umod{j+\ell}{r-\ell},\qquad
 \Om^2\Umod{j}{\ell}\cong \Umod{j+r}{\ell}.
\]
The least syzygy period of \(\Umod{j}{\ell}\) is
\[
 \operatorname{per}_{\Om}\Umod{j}{\ell}=
 \begin{cases}
 n/\gcd(n,\ell), & r=2\ell,\\
 2n/\gcd(n,r), & r\ne2\ell.
 \end{cases}
\]
The least positive even return time is always \(2n/\gcd(n,r)\). The stable syzygy matrix is a permutation matrix.
\end{proposition}
\begin{proof}
The kernel of \(U(j,r)\twoheadrightarrow \Umod{j}{\ell}\) is its complementary tail, of length \(r-\ell\) and top \(j+\ell\). Applying this twice restores the length and advances the top by \(r\). If \(r\ne2\ell\), odd returns are impossible because they change the length. Even returns advance the top by a multiple of \(r\), giving period \(2n/\gcd(n,r)\). If \(r=2\ell\), every syzygy preserves length and advances the top by \(\ell\), so the least period is \(n/\gcd(n,\ell)\). Finally the map \((j,\ell)\mapsto(j+\ell,r-\ell)\) is a permutation of the finite state set.
\end{proof}

\begin{example}[even and odd periods]\label{ex:orbit}
For \((n,r)=(5,3)\), the orbit of \(S_0\) begins
\[
 U(0,1)\longmapsto U(1,2)\longmapsto U(3,1)
 \longmapsto \Umod{4}{2}\longmapsto U(1,1).
\]
Its length is ten: every second step advances the top by three, and an odd return cannot restore the length. For \((n,r)=(6,4)\), the module \(\Umod{0}{2}\) instead has period three:
\[
 \Umod{0}{2}\xrightarrow{\Om}\Umod{2}{2}
 \xrightarrow{\Om}\Umod{4}{2}
 \xrightarrow{\Om}\Umod{0}{2}.
\]
Its least positive even return time is six.
\end{example}

\subsection{Maps and extensions from eligible depths}

A map between two uniserial modules is determined by the image of a top generator. In the target, the allowed path-basis vectors give a basis of maps. We number their depths from zero at the top. The vertex label gives a congruence, and the annihilation condition gives an interval of allowed depths. To compute Ext, we remove the maps that extend to the preceding projective in the resolution.

\begin{theorem}[extensions over a truncated cycle]\label{thm:all-cyclic-ext}
Let \(X=\Umod{j}{\ell}\) and \(Y=U(k,m)\), where \(1\le\ell<r\) and \(1\le m\le r\). Write \(\Om^iX=U(a_i,b_i)\); explicitly,
\[
 (a_{2t},b_{2t})=(j+tr,\ell),\qquad
 (a_{2t+1},b_{2t+1})=(j+tr+\ell,r-\ell).
\]
For every \(i\ge1\),
\[
 \kdim\Ext_{\Lambda_{n,r}}^i(X,Y)
 =\#\left\{q\in\mathbb Z:
 \begin{array}{l}
 \max(0,m-b_i)\le q<\min(m,r-b_i),\\
 q\equiv a_i-k\pmod n
 \end{array}\right\}.
\]
For arbitrary \(1\le\ell,m\le r\), degree zero is given by
\[
\begin{aligned}
 \kdim\Hom(\Umod{j}{\ell},U(k,m))=\#\{q\in\mathbb Z:\,&\max(0,m-\ell)\le q<m,\\
 &q\equiv j-k\pmod n\}.
\end{aligned}
\]
For a projective source, positive-degree extensions vanish.
\end{theorem}
\begin{proof}
Let \(y_0,\ldots,y_{m-1}\) be the usual path-basis vectors of \(Y\), with \(y_q\) at depth \(q\) and vertex \(k+q\). A homomorphism from \(U(a,b)\) is determined by the image of its top generator. Such an image must lie at vertex \(a\), and it must be annihilated by the path of length \(b\). Thus a basis for the possible images is given by
\[
 y_q\quad\text{with}\quad q\equiv a-k\pmod n,\qquad \max(0,m-b)\le q<m.
\]
These vectors give independent homomorphisms, even when the same vertex occurs several times in \(Y\). Taking \(a=j\) and \(b=\ell\) proves the Hom formula.

To calculate extensions, put \(d_i=r-b_i\). The embedding of \(\Om^iX=U(a_i,b_i)\) into the preceding projective \(P_{a_i-d_i}\) begins at depth \(d_i\). A map from that projective is determined by an arbitrary vector of \(Y\) at vertex \(a_i-d_i\). On a basis vector \(y_s\), restriction multiplies by the path of length \(d_i\), giving \(y_{s+d_i}\) if \(s+d_i<m\), and zero otherwise. Thus the image of restriction is spanned by exactly the eligible basis maps with \(q\ge d_i\). Conversely, each such map lifts by sending the projective's top generator to \(y_{q-d_i}\). The cokernel in \eqref{eq:cokernel} therefore retains precisely the eligible depths \(q<d_i\). Combining this with the annihilation bound and applying dimension shift gives the formula.
\end{proof}

The dimensions can exceed one: over \(\kk C_2/J^7\),
\[
 \kdim\Ext^1(\Umod{0}{3},U(1,3))=2.
\]
Here \(\Om U(0,3)=U(1,4)\), whose generator can map to either circled occurrence in the target:
\[
\begin{tikzpicture}[baseline=(current bounding box.center)]
\node[draw,circle,inner sep=3pt] (y0) at (0,0) {\(1\)};
\node (y1) at (2,0) {\(0\)};
\node[draw,circle,inner sep=3pt] (y2) at (4,0) {\(1\)};
\draw[->] (y0)--(y1); \draw[->] (y1)--(y2);
\node at (0,-.65) {\(q=0\)}; \node at (2,-.65) {\(q=1\)}; \node at (4,-.65) {\(q=2\)};
\node at (-1.7,0) {\(U(1,3):\)};
\end{tikzpicture}
\]
Both maps survive in Ext: the preceding projective contains the syzygy at depth three, so its restriction to this length-three target is zero.

Every positive-degree sequence has period dividing \(2n/\gcd(n,r)\), but its least period can be shorter than the source's syzygy period. We now determine the least periods for the arrow-ideal pairs considered in the dissertation.

%

\subsection{Least periods for arrow ideals}

\begin{corollary}[arrow-ideal dimensions and least periods]\label{thm:cyclic-ext}
Assume \(n\ge3\) and \(2\le r\le n\), and write \(R(a_j)=U(j+1,r-1)\). Put \(\varphi(i)=\lfloor i/2\rfloor r+(i\bmod2)\). For \(i\ge1\),
\[
 \kdim\Ext_{\Lambda_{n,r}}^i(R(a_j),R(a_l))
   =[l-j\equiv\varphi(i)\pmod n].
\]
Let \(g=\gcd(n,r)\). If \(r=2\), the least period is \(n\). If \(r\ge3\) and \(l-j\) is congruent to \(0\) or \(1\) modulo \(g\), the least period is \(2n/g\). Otherwise all positive-degree terms are zero. In degree zero,
\[
 \kdim\Hom(R(a_j),R(a_l))
   =[(j-l)\bmod n\in\{0,\ldots,r-2\}].
\]
\end{corollary}
\begin{proof}
\propref{prop:cyclic} gives
\[
 \Om^{2m}R(a_j)=R(a_{j+mr}),\qquad
 \Om^{2m+1}R(a_j)=S_{j+(m+1)r}.
\]
Substituting into \thmref{thm:all-cyclic-ext} leaves one possible depth on each parity and gives the displayed congruence. The even degrees attain exactly the class \(0\) modulo \(g\), and the odd degrees attain class \(1\).

On each active parity, the indicator has least period \(n/g\) in the half-degree. Hence the smallest even period is \(2n/g\). An odd period \(2t+1\) would require both parities to be active. Shifting a nonzero even term by this period changes \(\varphi\) by \(tr+1\); shifting a nonzero odd term changes it by \((t+1)r-1\). Since both shifted terms must still be nonzero, we obtain
\[
tr+1\equiv0,\qquad (t+1)r-1\equiv0\pmod n.
\]
Subtracting forces \(r\equiv2\pmod n\), hence \(r=2\) under \(2\le r\le n\). In that case \(\varphi(i)=i\). The Hom formula follows from the same theorem in degree zero.
\end{proof}

\begin{example}[reading the congruence]
Over \(\Lambda_{5,3}\), the self-pair \(R(a_0),R(a_0)\) has nonzero positive extensions exactly in degrees \(0\) and \(7\) modulo ten. The pair \(R(a_0),R(a_1)\) has them in degrees \(1\) and \(4\) modulo ten:
\[
\begin{array}{c|rrrrrrrrrr}
i&1&2&3&4&5&6&7&8&9&10\\ \hline
\dim\Ext^i(R(a_0),R(a_0))&0&0&0&0&0&0&1&0&0&1\\
\dim\Ext^i(R(a_0),R(a_1))&1&0&0&1&0&0&0&0&0&0
\end{array}
\]
For \(\Lambda_{6,3}\), \(g=3\). The pair \(R(a_0),R(a_2)\) has no positive extensions because \(2\) is neither \(0\) nor \(1\) modulo \(3\). The source still has periodic syzygies: the target detects none of the relevant first extensions.
\end{example}

For \(r\ge3\) these algebras are not gentle, but their one-arm covers still force periodic syzygies. For split Nakayama algebras, Erdmann, Klász and Marczinzik show that a nonzero self-extension in degree one of an indecomposable module forces nonzero self-extensions in every positive degree \cite{EKM23}. In the restricted range of \corref{thm:cyclic-ext}, \(\Ext^1(R(a_j),R(a_j))=0\). Outside that range, take \(\kk C_3/J^6\) and \(X=\Umod{0}{3}\). Here \(\Om X\cong X\), and \thmref{thm:all-cyclic-ext} gives \(\kdim\Ext_A^i(X,X)=1\) for every \(i\ge1\), in agreement with the cited result.

\section{Polynomial growth from linked cycles}\label{sec:linked}

The growth mechanism is already visible in a radical-square-zero algebra. Take the quiver
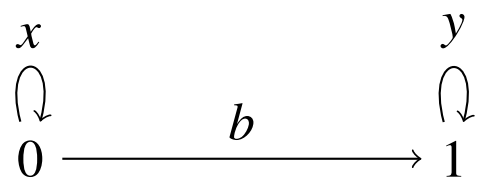
\[
\begin{tikzpicture}[baseline=-.4ex]
\node (u0) at (0,0) {\(0\)};
\node (u1) at (2.8,0) {\(1\)};
\draw[->] (u0) to[loop above] node {\(x\)} (u0);
\draw[->] (u1) to[loop above] node {\(y\)} (u1);
\draw[->] (u0)--(u1) node[midway,above] {\(b\)};
\end{tikzpicture}
\qquad A=\kk Q/J^2.
\]
The loops keep the simples where they are, while the arrow \(b\) produces a new summand:
\[
\Om S_0\cong S_0\oplus S_1,\qquad \Om S_1\cong S_1.
\]
Hence \(\Om^iS_0\cong S_0\oplus S_1^{\oplus i}\), and for \(T=S_0\oplus S_1\),
\[
\kdim\Ext_A^i(S_0,T)=i+1\qquad (i\ge0).
\]
This is a two-vertex version of \exref{ex:polynomial}. The same mechanism works when the persistent summands are longer uniserial modules rather than simples.

\subsection{Two cycles}
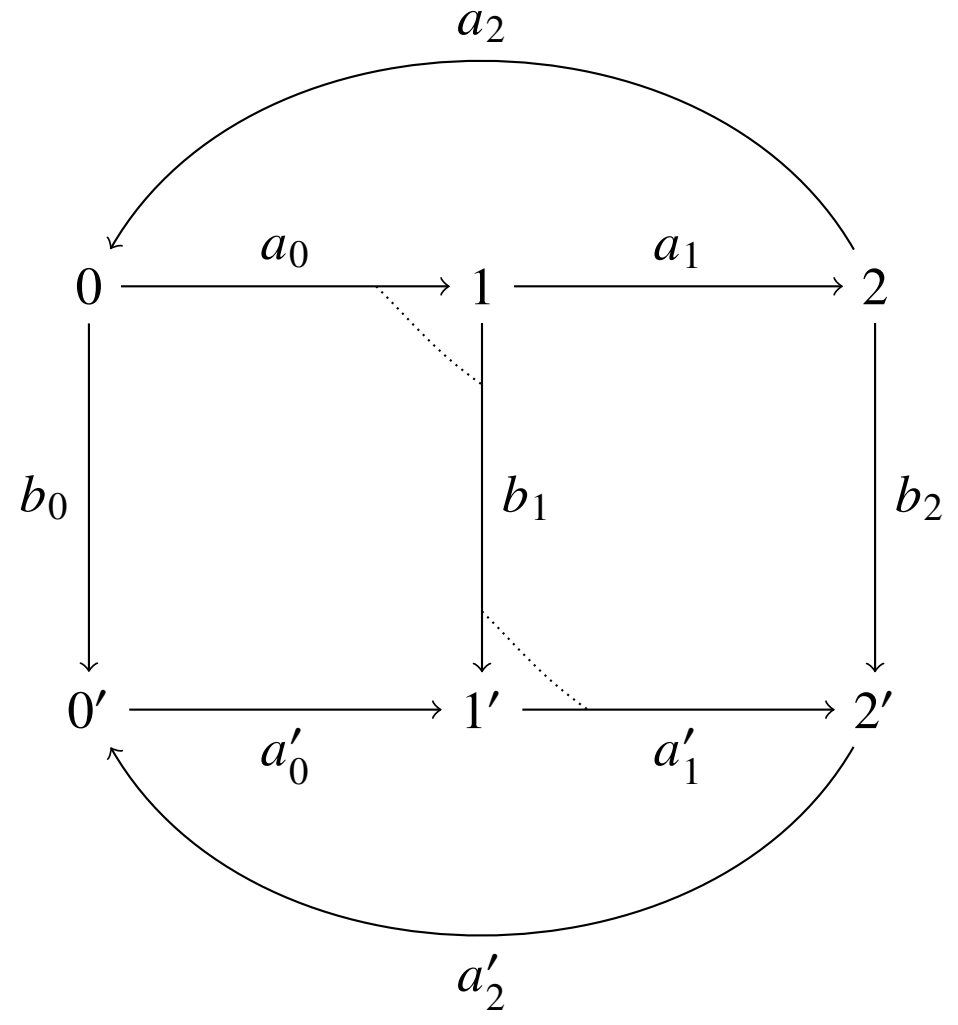
\begin{figure}[ht]
\centering
\begin{tikzpicture}
\node (u0) at (0,1.6) {\(0\)};
\node (u1) at (2.6,1.6) {\(1\)};
\node (u2) at (5.2,1.6) {\(2\)};
\node (v0) at (0,-1.2) {\(0'\)};
\node (v1) at (2.6,-1.2) {\(1'\)};
\node (v2) at (5.2,-1.2) {\(2'\)};
\draw[->] (u0)--(u1) node[midway,above] {\(a_0\)};
\draw[->] (u1)--(u2) node[midway,above] {\(a_1\)};
\draw[->] (u2) to[out=120,in=60] node[above] {\(a_2\)} (u0);
\draw[->] (v0)--(v1) node[midway,below] {\(a'_0\)};
\draw[->] (v1)--(v2) node[midway,below] {\(a'_1\)};
\draw[->] (v2) to[out=-120,in=-60] node[below] {\(a'_2\)} (v0);
\draw[->] (u0)--(v0) node[midway,left] {\(b_0\)};
\draw[->] (u1)--(v1) node[midway,right] {\(b_1\)};
\draw[->] (u2)--(v2) node[midway,right] {\(b_2\)};
\draw[densely dotted] (1.9,1.6) to[out=-45,in=145] (2.6,.95);
\draw[densely dotted] (2.6,-.55) to[out=-45,in=145] (3.3,-1.2);
\end{tikzpicture}
\caption{Two linked \(3\)-cycles; the dotted arcs mark \(a_0b_1=b_1a'_1=0\).}
\label{fig:doubled}
\end{figure}

Let \(C_n\) be an oriented cycle with vertex set \(\mathbb Z/n\mathbb Z\). We take two copies, with vertices \(p\) and \(p'\) and arrows
\[
a_p:p\to p+1,\qquad a'_p:p'\to(p+1)',\qquad b_p:p\to p'\quad(p\in\mathbb Z/n\mathbb Z).
\]
All indices are modulo \(n\); \figref{fig:doubled} shows the complete quiver for \(n=3\).

Kill all paths of length \(r\) along the upper cycle and of length \(s\) along the lower cycle, where \(3\le r,s\le n\). In addition impose, for every \(p\),
\[
a_{p-1}b_p=0,\qquad b_pa'_p=0.
\]
We denote the resulting string algebra by \(\Lambda_{n,r,s}\). The mixed relations ensure that \(b_pA\cong S_{p'}\).

Write \(U(p,\ell)\) for the uniserial module of length \(1\le\ell<r\) with top \(S_p\), supported on the upper cycle, and \(U(p',\ell)\) for its lower counterpart, with \(1\le\ell<s\). All links act as zero. The upper projective cover has an extra arm \(b_p\), so
\[
\begin{aligned}
 \Om U(p,\ell)&\cong U(p+\ell,r-\ell)\oplus S_{p'},\\
 \Om U(p',\ell)&\cong U((p+\ell)',s-\ell).
\end{aligned}
\]
Retaining only the complementary tail permutes the modules on each row. On this family the syzygy matrix therefore has the form
\[
 B_{\rm cyc}=\begin{pmatrix}P&C\\0&Q\end{pmatrix},\qquad
(B_{\rm cyc}^m)_{12}=\sum_{a=0}^{m-1}P^aCQ^{m-1-a}\quad(m\ge1),
\]
where \(P,Q\) are permutation matrices and \(C\) records the link simples. The orbits of \(P\) and \(Q\) are the recurrent components of the syzygy graph. The summation counts the possible times of transfer from an upper orbit to a lower one.

After two syzygies these are all the modules we need. By \propref{prop:reduction}, the summands are positive-length path ideals: a path along either cycle leaves a tail shorter than its truncation length, and a link gives a lower simple.

A walk meets at most two orbits, one on each row. By \thmref{thm:growth}, every extension sequence is eventually zero or eventually quasi-polynomial of degree at most one. If all summands of \(\underline{\Om^2M}\) lie on the lower cycle, the sequence is bounded. A linear term occurs on some residue class precisely when an upper state in \(\operatorname{supp}v_M\) can reach a lower state in \(\operatorname{supp}h_N\).

\begin{example}[one new summand at every step]\label{ex:doubled-count}
Take \(n=r=s=3\). The first syzygies of the upper simple \(S_0\) are
\[
\begin{aligned}
\Om S_0&\cong U(1,2)\oplus S_{0'},\\
\Om^2S_0&\cong S_0\oplus S_{1'}\oplus U(1',2),\\
\Om^3S_0&\cong U(1,2)\oplus S_{0'}^{\oplus2}\oplus U(2',2).
\end{aligned}
\]
At each step the upper summand produces one new lower simple, while each existing lower summand has just one successor. Hence \(\Om^iS_0\) has \(i+1\) indecomposable summands, all with simple top. Taking \(T=A/\rad A\), \lemref{lem:homological} gives
\[
 \kdim\Ext_A^i(S_0,T)=i+1\qquad(i\ge0).
\]
\end{example}

\subsection{Several cycles and accessible links}
Take finitely many disjoint oriented cycles \(C_u\), with \(n_u\) vertices, and kill their paths of length \(r_u\), where \(3\le r_u\le n_u\). Join them by arrows, called links, pointing forward in some ordering of the cycles. Each vertex may be the source of at most one link and the target of at most one link. Kill every length-two path containing a link. The resulting algebra is string, and each link \(b:c\to d\) gives \(bA\cong S_d\).

Write \(U_u(j,\ell)\) for the uniserial module of length \(1\le\ell<r_u\) with top \(S_j\) on \(C_u\), with all links acting as zero. Its cover gives
\begin{equation}\label{eq:linked-cover}
 \Om U_u(j,\ell)\cong U_u(j+\ell,r_u-\ell)
       \oplus\!\bigoplus_{b:j\to d\ {\rm a\ link}}S_d.
\end{equation}
All summands are non-projective. Retaining the summand on \(C_u\) gives the permutation
\[
 F_u(j,\ell)=(j+\ell,r_u-\ell),\qquad F_u^2(j,\ell)=(j+r_u,\ell),
\]
with vertex indices modulo \(n_u\). Starting from \(S_w\), its successive tops are
\[
w,\ w+1,\ w+r_u,\ w+r_u+1,\ w+2r_u,\ldots.
\]
Thus, if \(g_u=\gcd(n_u,r_u)\), the source \(c\) of a link occurs exactly when
\begin{equation}\label{eq:live}
 c\equiv w\quad\hbox{or}\quad c\equiv w+1\pmod{g_u}.
\end{equation}
Every such top recurs periodically. Using that link starts a new persistent summand at its target vertex on the next cycle.

\begin{definition}[live chains]\label{def:live}
A chain of links is \emph{live from \(v\)} if \eqref{eq:live} holds successively, beginning at \(v\) and then at each preceding link target. Write \(\delta^*(v)\) for the largest number of links in such a chain, with value zero if none can be used.
\end{definition}

\begin{proposition}[the degree from a live chain]\label{prop:linked}
Let \(s_i(v)\) count the indecomposable summands of \(\underline{\Om^iS_v}\), with multiplicity. Then \(s_i(v)\) is eventually quasi-polynomial of degree \(\delta^*(v)\).
\end{proposition}
\begin{proof}
By \eqref{eq:linked-cover}, syzygies stay within direct sums of the finitely many modules \(U_u(j,\ell)\). Since links point forward, the recurrent components of the syzygy graph on these modules are exactly the orbits of the permutations \(F_u\). A transfer between reachable orbits uses a live link and enters the next orbit at a simple. The greatest number of orbits met by a walk starting at \(S_v\) is therefore \(\delta^*(v)+1\). Apply \thmref{thm:growth} with terminal weight \(\onevec\).
\end{proof}

For a pair \((M,N)\), start instead at the summands of \(\underline{\Om^2M}\). From \(U_u(j,\ell)\), the first link is accessible when its source is congruent to \(j\) or \(j+\ell\) modulo \(g_u\); all subsequent transfers start from simples and use \eqref{eq:live}. Finally discard states unable to reach \(\operatorname{supp}h_N\). The last restriction can lower the degree: a target need not detect every summand.

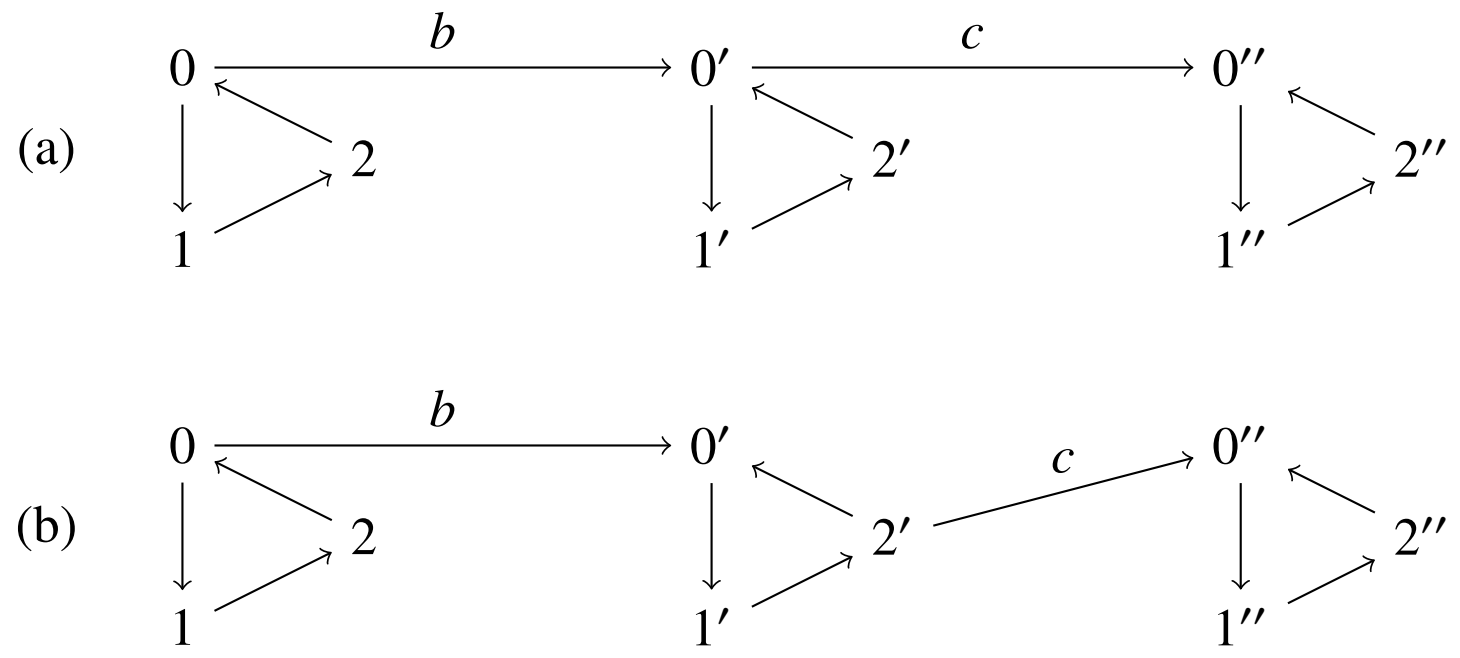
\begin{figure}[ht]
\centering
\begin{tikzpicture}
\begin{scope}
\node at (-.9,-.55) {(a)};
\node (a0) at (0,0) {\(0\)}; \node (a1) at (0,-1.2) {\(1\)}; \node (a2) at (1.2,-.6) {\(2\)};
\node (b0) at (3.5,0) {\(0'\)}; \node (b1) at (3.5,-1.2) {\(1'\)}; \node (b2) at (4.7,-.6) {\(2'\)};
\node (c0) at (7,0) {\(0''\)}; \node (c1) at (7,-1.2) {\(1''\)}; \node (c2) at (8.2,-.6) {\(2''\)};
\draw[->] (a0)--(a1); \draw[->] (a1)--(a2); \draw[->] (a2)--(a0);
\draw[->] (b0)--(b1); \draw[->] (b1)--(b2); \draw[->] (b2)--(b0);
\draw[->] (c0)--(c1); \draw[->] (c1)--(c2); \draw[->] (c2)--(c0);
\draw[->] (a0)--(b0) node[midway,above] {\(b\)};
\draw[->] (b0)--(c0) node[midway,above] {\(c\)};
\end{scope}
\begin{scope}[shift={(0,-2.5)}]
\node at (-.9,-.55) {(b)};
\node (a0) at (0,0) {\(0\)}; \node (a1) at (0,-1.2) {\(1\)}; \node (a2) at (1.2,-.6) {\(2\)};
\node (b0) at (3.5,0) {\(0'\)}; \node (b1) at (3.5,-1.2) {\(1'\)}; \node (b2) at (4.7,-.6) {\(2'\)};
\node (c0) at (7,0) {\(0''\)}; \node (c1) at (7,-1.2) {\(1''\)}; \node (c2) at (8.2,-.6) {\(2''\)};
\draw[->] (a0)--(a1); \draw[->] (a1)--(a2); \draw[->] (a2)--(a0);
\draw[->] (b0)--(b1); \draw[->] (b1)--(b2); \draw[->] (b2)--(b0);
\draw[->] (c0)--(c1); \draw[->] (c1)--(c2); \draw[->] (c2)--(c0);
\draw[->] (a0)--(b0) node[midway,above] {\(b\)};
\draw[->] (b2)--(c0) node[midway,above] {\(c\)};
\end{scope}
\end{tikzpicture}
\caption{Growth from \(S_0\): (a) quadratic; (b) linear.}
\label{fig:live}
\end{figure}

\begin{example}[a live link and a missed link]\label{ex:live}
Use the quivers in \figref{fig:live}, killing all paths of length three along a cycle and all paths of length two containing a link. Starting at \(S_0\), the first link produces \(S_{0'}\). Following only the summand that stays on the middle cycle gives
\[
S_{0'}\longmapsto U(1',2)\longmapsto S_{0'}.
\]
The tops are \(0'\) and \(1'\). The second link in (a) is therefore used repeatedly, while the link in (b), starting at \(2'\), is never used. Although \(U(1',2)\) also has a composition factor at \(2'\), only its top produces a link summand in \eqref{eq:linked-cover}. The two summand counts therefore have degrees two and one, respectively.

The first few counts make the difference visible:
\[
\begin{array}{c|rrrrrrrrr}
i&0&1&2&3&4&5&6&7&8\\ \hline
\text{second link at }0'&1&2&3&4&6&7&10&11&15\\
\text{second link at }2'&1&2&2&3&3&4&4&5&5
\end{array}
\]
For the first row, in degree \(2m\) the three cycles contribute \(1\), \(m\), and \(\binom{m+1}{2}\) summands. In degree \(2m+1\), the middle count increases to \(m+1\) and the other two stay the same. In the second row, the third cycle contributes nothing and the second contributes \(\lceil i/2\rceil\) in degree \(i\). These numbers also give \(\kdim\Ext_A^i(S_0,A/\rad A)\).
\end{example}

\printbibliography

\section*{Acknowledgements}
I thank  Azzurra Ciliberti and Karin Baur for their guidance and supervision during my dissertation.

\end{document}